# GPU ACCELERATION OF THE SEVEN-LEAGUE SCHEME FOR LARGE TIME STEP SIMULATIONS OF STOCHASTIC DIFFERENTIAL EQUATIONS

SHUAIQIANG LIU, GRAZIANA COLONNA, LECH GRZELAK AND CORNELIS OOSTERLEE

**Abstract** Monte Carlo simulation is widely used to numerically solve stochastic differential equations. Although the method is flexible and easy to implement, it may be slow to converge. Moreover, an inaccurate solution will result when using large time steps. The Seven League scheme [1], a deep learning-based numerical method, has been proposed to address these issues. This paper generalizes the scheme regarding parallel computing, particularly on Graphics Processing Units (GPUs), improving the computational speed.

## Introduction

In this paper, we will develop a highly accurate numerical discretization scheme for stochastic differential equations (SDEs), which is based on taking possibly large discrete time steps. We "learn" to take large time steps [1] by using an artificial neural network (ANN), in the context of supervised machine learning, with the help of stochastic collocation polynomials (SCMC, see [2]).

In many forms and flavors, the deep learning paradigm [3] receives much attention in science and engineering nowadays. The Physics-Informed Neural Networks (PiNN) [4], for example, combining physical and mathematical insights in an unsupervised learning fashion, starts to enter the field of solving ordinary (ODEs) and partial differential equations (PDEs). The corresponding computational costs are nontrivial, however, simply because the underlying equations need to be learned from scratch, and this costs time. Supervised learning, based on labels, on the other hand, is a classical form of machine learning, which is often more efficient as there is an offline stage, in which the input-output labelled relations are being learned, followed by a highly efficient online stage, where the learned manifold of solutions is evaluated for new input values. In our work, supervised learning is employed, which, however, strongly relies on a sophisticated numerical stochastic collocation technique, to achieve a challenging task within numerical analysis. In the present paper, we want to show that by parallelization the method's efficiency can be greatly enhanced, thus computing times of both the online and the offline stages are significantly reduced, on a graphics processing unit (GPU).

The basic idea of the 7L scheme is to learn a small number of (conditional) stochastic collocation points, and the learned neural network function is employed to forecast the unknown collocation points for the next time step. Then, by means of the stochastic collocation Monte Carlo sampler (SCMC) [2], the random paths are generated. Interestingly, the strong convergence error of the 7L scheme is independent of the size of the simulation time step. In other words, different from classical numerical schemes, the accuracy of the numerical solution does not decrease when solving SDEs with large time steps due to the learning stage.

Parallel computing on GPUs is prevalent nowadays, especially for numerical simulations that require intensive computational resources. Regarding stochastic differential equations, parallel implementations of classical Monte Carlo simulations on GPUs have been well studied, for example, [5, 6]. In this paper, we extend the 7L scheme by parallel computing to further improve its computational speed. There is inherent parallelism in the 7L scheme, as a neural network consists of a large number of artificial neurons that can work in parallel. Moreover, required interpolation functions are independent of each other. Thus, they can be easily distributed over different processing units.

## Methodology, the 7L Scheme

Suppose a real-valued random variable $Y(t)$ is defined on the probability space $(\Omega, \Sigma, \mathbb{P})$ with filtration $\mathcal{F}_{t \in [0,T]}$, sample space $\Omega$, $\sigma$-algebra $\Sigma$ and probability measure $\mathbb{P}$. For the time evolution of $Y(t)$, consider the generic scalar Itô SDE,

$$dY(t) = a(t, Y(t), \theta)dt + b(t, Y(t), \theta)dW(t), \quad 0 \leq t \leq T, \quad (6.1)$$

with the drift term $a(t, Y(t), \theta)$, the diffusion term $b(t, Y(t), \theta)$, model parameters $\theta$, Wiener process $W(t)$, and given initial value $Y_0 := Y(t=0)$. The solution of (6.1) is unique when the drift and diffusion terms meet some regularity conditions.

The basic discretization for each Monte Carlo path, is the Euler-Maruyama scheme [7], which reads,

$$\hat{Y}_{i+1} | \hat{Y}_i = \hat{Y}_i + a(t_i, \hat{Y}_i, \theta)\Delta t + b(t_i, \hat{Y}_i, \theta)\sqrt{\Delta t}\hat{X}_{i+1}, \quad (6.2)$$

where $\hat{Y}_{i+1} := \hat{Y}(t_{i+1})$ is a realization (i.e., a number) from random variable $\tilde{Y}(t_{i+1})$, which represents the numerical approximation to exact solution $Y(t_{i+1})$ at time point $t_{i+1}$, and a realization $\hat{X}_{i+1}$ is drawn from the random variable $X$, which here follows the standard normal distribution $\mathcal{N}(0, 1)$. The Euler-Maruyama scheme will be used to generate the training data set. However, training will be based on tiny time steps (for accuracy reasons).

Similarly, the 7L scheme reads $\hat{Y}_{i+1} | \hat{Y}_i = g_m(\hat{X}_{i+1})$, where $g_m(\cdot)$ stands for a mapping function, transforming a standard normal distribution to the target distribution at time $t_{i+1}$. The function $g_m(\cdot)$ can be obtained through an interpolation technique, based on $m$ pairs of collocation points $(x_j, y_j)$, where $j = 0, \ldots, m-1$, $x_j$ are obtained from the standard norm distribution $X$ (here Gauss-Hermite quadrature points), and $y_j$ are stochastic collocation points at time $t_{i+1}$, *conditional* on the previous realization $\hat{Y}_i$. In the context of Markov processes, the function of computing $y_j$ can be written as follows,

$$y_j(t_{i+1}) | \hat{Y}_i = H_j\left(\hat{Y}_i, t_{i+1} - t_i, \theta\right). \quad (6.3)$$

So, a neural network can be trained to approximate the function $H_j(\cdot)$. The two key components of the scheme refer to the function $H_j(\cdot)$ and the interpolation function $g_m(\cdot)$, both of which will be parallelized in Section 6.2.

## Parallelization

Parallelization is carried out by the parallel implementation of the appearing interpolation functions and the neural network involved in the algorithm.

Note that the parallelization of the 7L scheme is focused on the online stage, because the training stage is done offline and once. A variant, the 7L-CDC scheme, see the paper [1], with more interpolations in Step 3 above, can be parallelized similarly. In this work, we use the barycentric version [8] of Lagrange interpolation on GPUs and CPUs to fairly compare their speed performance.

Please refer to the original paper [1] for more interpolation techniques.

**Algorithm I: A parallel algorithm of the 7L scheme**

1. **Offline stage**: Train the ANNs to predict the stochastic collocation points. At this stage, we choose different $\theta$ values, simulate corresponding Monte Carlo paths, with small constant time increments $\Delta\tau = \tau_{i+1} - \tau_i$ in $[0, \tau_{\max}]$, generate the corresponding collocation points $\hat{y}_j|\hat{Y}_i \approx y_j|\hat{Y}_i$, and learn the relation between inputs and outputs to obtain $\hat{H}_k \approx H_k$.
2. **Online stage**: Partition time interval $[0,T]$, $t_i = iT/N$, $i = 0, \ldots N$, with equidistant "large" time step $\Delta t = t_{i+1} - t_i$, and output $N$ sample paths.
3. Run the ANNs to compute $m$ collocation points at time $t_{i+1}$ for each path,

$$\hat{y}_j(t_{i+1})|\hat{Y}_i = \hat{H}_j(\hat{Y}_i, t_{i+1} - t_i, \theta), j = 1, 2, \ldots, m, \quad (6.4)$$

   and form a vector
   $\hat{y}_{i+1} = (\hat{y}_1(t_{i+1})|\hat{Y}_i, \hat{y}_2(t_{i+1})|\hat{Y}_i, \ldots, \hat{y}_m(t_{i+1})|\hat{Y}_i)$.
   This step is parallelized by running the ANNs in a batch model on the GPU.
4. Divide $N$ sample paths into $N_T$ groups, and allocate a group of $\frac{N}{N_T}$ paths to a certain thread on GPUs.
5. For each of $\frac{N}{N_T}$ paths in a thread, compute interpolation function $g_m(\cdot)$, based on $m$ pairs of $(x_j, \hat{y}_j)$.
6. Sample from $X$ and obtain a sample $\hat{Y}_{i+1}|\hat{Y}_i = g_m(\hat{X}_{i+1})$.
7. Collect all paths $\hat{Y}_{i+1}$ from $N_T$ threads to form a compete set at time $t_{i+1}$.
8. Return to Step 3 by $t_{i+1} \to t_i$, iterate until terminal time $T$.

## Numerical Results

In this section, we evaluate the computational performance of the parallized 7L scheme. Here we take the Ornstein-Uhlenbeck (OU) process as an example. The OU process is a mean reverting process, defined as follows,

$$dY(t) = -\lambda(Y(t) - \overline{Y})dt + \sigma dW(t), \quad 0 \leq t \leq T, \quad (6.5)$$

with $\overline{Y}$ the long term mean of $Y(t)$, $\lambda$ the speed of mean reversion, and $\sigma$ the volatility. The initial value is $Y_0$, and the model parameters are $\theta := \{\overline{Y}, \sigma, \lambda\}$. Its analytical solution is given by,

$$Y(t) \stackrel{d}{=} Y_0 e^{-\lambda t} + \overline{Y}(1 - e^{-\lambda t}) + \sigma\sqrt{\frac{1 - e^{-2\lambda t}}{2\lambda}}X, \quad (6.6)$$

with $t_0 = 0$, $X \sim \mathcal{N}(0, 1)$. Equation (6.6) is used to compute the reference value to the path-wise error and the strong convergence.

In the training phase, the Euler-Maruyama scheme is used to discretize the OU dynamics and generate the data set (here five stochastic collocation points to learn within the ANN). After the training, the 7L scheme with the obtained ANNs is used to solve the OU process, as shown in **Algorithm I**.

The ANN hyperparameters are set as follows here. We use 4 hidden layers, 50 neurons per layers, a Softplus activation function, a Glorot initialization, the Adam optimizer, a batch size of 1024, and a learning rate of $10^{-3}$.

The parallized 7L scheme is evaluated on the GPU and CPU as follows,

GPU  Type: GeForce MX150, Graphic Cores: 384, Graphics clock: 1468 MHz, Memory speed: 6.01 GHz, Memory bandwidth: 48.06 Gb/s.

CPU  Type: Intel Core i7-8550U, Cores: 4, Maximum speed: 4.0 GHz, Base clock speed: 1.80 GHz.

The parallelization is done in CUDA (Compute Unified Device Architecture), the platform created by Nvidia. The threads, which are the basic operational units in CUDA, are computing processes that run in parallel. The number of threads $N_T$ used for the parallelization is 256. Therefore, the number of paths $N_B$ per thread is proportional to the total number of simulated paths $N_P$, $N_B = \frac{N_P}{N_T}$. The speedup ratio is defined by the running time of the parallelized code (on the GPU) divided by the running time of the original code (on the CPU). The reported time is obtained by running the corresponding code 100 times and taking the averaged execution time.

Table 6.1: Computational time (seconds) of the 7L and 7L-CDC schemes.

| Number of paths | the 7L scheme | | | the 7L-CDC scheme | | |
| --- | --- | --- | --- | --- | --- | --- |
| | Sequential time | Parallel time | Speedup | Sequential time | Parallel time | Speedup |
| 1,000 | 1.555 | 0.268 | 5.8 | 1.296 | 0.062 | 20.9 |
| 50,000 | 70.108 | 6.844 | 10.2 | 64.731 | 2.828 | 22.9 |
| 100,000 | 134.745 | 14.623 | 9.2 | 132.198 | 5.886 | 22.5 |
| 200,000 | 282.456 | 21.545 | 13.1 | 251.684 | 11.527 | 21.8 |

As shown in Table 6.1, the speedup ratio appears to converge to 10 when the number of the sample paths increases. However, this ratio may fluctuate due to the unstable performance of GPUs and CPUs. The major acceleration comes from the parallelization of the interpolation process (here based on five collocation points), as the ANN running times on the GPU and CPU have a small difference in this test. The 7L-CDC scheme employs a global interpolation technique which is based on the marginal collocation points to compute the conditional collocation points for each random path, instead of using the ANNs for each path as the 7L scheme does in Step 3 of **Algorithm I**. The 7L-CDC scheme only requires the ANNs to compute a small number of marginal collocation points (here five marginal collocation points) along with the probability distribution. The 7L-CDC scheme can be used as a faster variant of the 7L scheme, as long as the global interpolation technique is computationally cheaper than the evaluation of the ANN. As shown in Table 6.1, the speedup ratio of the parallelized 7L-CDC scheme converges to 22. There are two interpolation processes (one for the five conditional collocation points and another for five marginal collocation points) in the 7L-CDC scheme, which explains why the speedup ratio is as twice as that of the 7L scheme.

The speedup ratio is also affected by other factors, for example, the number of threads and the configuration of the used GPU. The original paper [1] has proved that the numerical error does not grow when the simulation time step size increases. We find that the above property holds when the 7L scheme is implemented on GPUs in a parallel way.

Summarizing, a neural networks-based numerical solver for stochastic differential equations, the 7L scheme, has been parallelized to accurately carry out large time step simulations, with a further computational acceleration by a factor of 10 or even 20, on the used GPU.

**Acknowledgment**
This project is part of the ABC-EU-XVA project and has received funding from the European Unions Horizon 2020 research and innovation programme under the Marie Skldowska–Curie grant agreement No 813261.